\documentclass[12pt]{svmult}

\usepackage{amsfonts}
\usepackage{latexsym}
\usepackage{color}
 \RequirePackage{amsbsy} 
 \RequirePackage{amsopn} 
 \RequirePackage{amsmath}
  \RequirePackage{amssymb}
  \usepackage[mathscr]{eucal}

 \newtheorem{thee}{Theorem}
 \newtheorem{coor}[thee]{Corollary}
 \newtheorem{leem}[thee]{Lemma}
 
 \newtheorem{prro}[thee]{Proposition}

 \newtheorem{reem}[thee]{Remark}

 \newcommand{\balf}
 {\renewcommand{\theenumi}{(\alph{enumi})}
 \renewcommand{\labelenumi}{\theenumi}
                      \begin{enumerate}}
\newcommand{\ealf}   {\end{enumerate}
                      \renewcommand{\theenumi}{\arabic{enumi}}
                      \renewcommand{\labelenumi}{\theenumi.}}
\newcommand{\bara}   {\renewcommand{\theenumi}{(\arabic{enumi})}
                      \renewcommand{\labelenumi}{\theenumi}
                      \begin{enumerate} }
\newcommand{\eara}   {\end{enumerate}
                      \renewcommand{\theenumi}{\arabic{enumi}}
                      \renewcommand{\labelenumi}{\theenumi.}}

 \newcommand{\brom}   {\renewcommand{\theenumi}{(\roman{enumi})}
                      \renewcommand{\labelenumi}{\theenumi}
                      \begin{enumerate} }
\newcommand{\erom}   {\end{enumerate}
                      \renewcommand{\theenumi}{\arabic{enumi}}
                      \renewcommand{\labelenumi}{\theenumi.}}

	  \newcommand{\cl}{^{{\mbox{\rm \footnotesize{cl}}^{{\star}}}}}
	  \newcommand{\cls}{^{{\mbox{\rm 
	  \footnotesize{cl}}^{{\overline{\star}}}}}}
	  \newcommand{\Na}{\mbox{\rm Na}}
	  \newcommand{\Kr}{\mbox{\rm Kr}}
	   \newcommand{\KN}{\mbox{\rm KN}}
	   
	   \newcommand{\Max}{\mbox{\rm Max}}
	    \newcommand{\Pic}{\mbox{\rm Pic}}

  \DeclareMathOperator{\calM}     {\mathcal M}%
    \DeclareMathOperator{\calL} {\boldsymbol{\mathcal L}}%
\begin{document}

\title*{An historical overview of Kronecker function rings, Nagata 
rings, and related star and semistar operations}

\titlerunning{An historical overview of Kronecker function rings and related topics}

\author{Marco Fontana\inst{1}\and
K. Alan Loper\inst{2}}
\institute{Dipartimento di Matematica 
Universit\`a degli Studi Roma Tre \\
Largo San Leonardo Murialdo, 1 \\ 
00146 Roma, Italy
\texttt{fontana@mat.uniroma3.it}
\and Department of Mathematics \\ Ohio State University-Newark\\ Newark,
Ohio 43055 USA
 \texttt{lopera@math.ohio-state.edu}}
%
%
\maketitle

\section{Introduction: The Genesis }

 Toward the middle of the  XIXth century,  E.E. Kummer  
discovered that the 
 ring of integers of a cyclotomic field does not have the unique 
 factorization property and he introduced the concept of 
``ideal numbers'' to re-establish 
 some of the factorization theory for cyclotomic integers  
\cite[Vol. 1, 203-210, 583-629]{Kummer}.

 As  R. Dedekind   wrote in 1877 to his former student E. 
Selling,   \sl the 
 goal \rm  of \sl a \rm general theory was immediately clear 
 after Kummer's solution in the special case of cyclotomic integers: 
to 
 extend   Kummer's theory to the case of general algebraic integers.
 
 Dedekind admitted to having struggled unsuccessfully for many years 
before he published the first version of his theory in  1871 
\cite{Kummer}
(XI 
supplement to Dirichlet's ``Vorlesungen \"uber Zahlentheorie''  
\cite{Dedekind}).

The theory of Dedekind domains, as it is known today, is based on 
 Dedekind's original ideas and results. Dedekind's point of view is based on  
ideals (``ideal numbers'') for 
generalizing the algebraic numbers; he proved that, \sl   in 
the ring of the 
integers of an algebraic number field, each proper ideal factors 
uniquely into 
a product of prime ideals. \rm

 L. Kronecker essentially achieved  this goal in  
1859, but he published nothing until 
1882  \cite{Kronecker}.

Kronecker's theory holds in a larger context than that of rings 
of integers of algebraic numbers and solves a more general problem.   
The primary objective of his theory was to extend the set of elements 
and the concept of divisibility in such a way that any finite set of 
elements has a  GCD (greatest common divisor) in an extension of the 
original ring which still mirrors as closely as possible the ideal 
structure of the original ring.
It is probably for this reason that  the basic objects of  
Kronecker's  
theory  --corresponding to Dedekind's ``ideals''-- are called  
``divisors''.

Let $D_{_{\!0}}$ be a PID with quotient field $K_{_{\!0}}$ and let 
$K$ be a finite 
field extension of $K_{_{\!0}}$. Kronecker's \it divisors \rm 
 are essentially all 
the possible GCD's of finite sets of elements of  $K$ that are 
algebraic over $K_{_{\!0}}$; 
\  a divisor \rm   is  \it   integral \rm  if 
it is the GCD of a finite set of 
elements of the integral closure $D$ of $D_{_{\!0}}$ in $K$.

One of the key points of Kronecker's theory is that it is possible to 
give an explicit description of the ``divisors''.\  The divisors can 
be 
represented as equivalence classes of polynomials and a given 
polynomial in $D[X]$ represents the class of the integral  divisor 
associated 
with the set of its coefficients.

More precisely, we can give the following definition.
\smallskip

\noindent \bf The classical Kronecker function ring. \rm Let $D$ be 
as above.  
 \it The Kronecker function ring of $D$   \rm  is given by:
$$
\begin{array}{rl}
\Kr(D) := & \hskip -5pt  \left\{\frac{f}{g} \, \mid \,   f,g \in D[X] 
\,, \; g \neq 0\; \; \mbox{and} \; \; 
{\boldsymbol c}(f) \subseteq {\boldsymbol c}(g)\right\} \\ 
= & \hskip -5pt\left\{\frac{f'}{g'} \, \mid  \, f',g' \in D[X] \; \; 
\mbox{ and} \; \;
{\boldsymbol c}(g') = D\right\}\,,
\end{array} 
$$
(where $\,\boldsymbol{c}(h)\,$ denotes  \it  the   
content  \rm 
  of a
  polynomial $\,h \in D[X]$, i.e. the ideal of $D$ generated by the 
  coefficients of $h$).

  Note that we are assuming that 
  $D$ is a Dedekind domain (being the integral closure of 
$D_{_{\!0}}$, which is a PID, in a 
  finite field extension $K$ of the quotient field $K_{_{\!0}}$ of 
  $D_{_{\!0}}$ \cite[Theorem 41.1 and Theorem 37.8]{Gilmer:1972}). \\ 

In this case, for each polynomial $g \in D[X]$,  
${\boldsymbol c}(g)$ is an invertible ideal of $D$ and, by choosing a 
polynomial $u \in 
K[X]$ such that ${\boldsymbol c}(u) = ({\boldsymbol c}(g))^{-1}:= (D: 
{\boldsymbol c}(g))$, 
we have 
$f/g = uf/ug =f'/g'$, with $f' := uf,\, 
  g' := ug \in D[X]$ and thus ${\boldsymbol c}(g') = D$ (Gauss Lemma).

  \medskip

The fundamental properties of the Kronecker function ring are the 
following (cf. \cite[Chapter II]{Weyl}, \cite[Theorem 32.6 (for 
$\star$ equal to the identity star operation)]{Gilmer:1972}): 
\begin{enumerate}
	\bf \item[(1)] \sl 
 $\Kr(D)$ is a B\'ezout domain  (i.e. each finite set of elements, 
not all zero, 
has a GCD and the GCD can 
	be expressed as linear combination of these elements) and  $D[X] 
\subseteq \Kr(D) \subseteq 
	K(X)$ (in particular, the field of rational functions $K(X)$ is the 
quotient 
	field of $\Kr(D)$).
	
\bf \item[(2)] \sl Let $a_{_{\!0}}, a_{_{\!1}},\ldots\!, 
a_{_{\!n}}\in D$ and set $f := a_{_{\!0}} +a_{_{\!1}}X +\ldots\!+ 
a_{_{\!n}}X^n \in D[X]$, then:
$$
\begin{array}{rl}
		& \hskip -5pt (a_{_{\!0}}, a_{_{\!1}},\ldots\!, a_{_{\!n}})\Kr(D) = 
		f \Kr(D) \;\,  ({\mbox{thus, \rm GCD}}_{{\mbox{\footnotesize \rm 
Kr}}(D)} (a_{_{\!0}}, 
		a_{_{\!1}},\ldots \!, a_{_{\!n}}) \!= \! f)\,,   \\ 
		     
		& \mbox{$f$}\, \! \Kr(D) \cap K = (a_{_{\!0}}, 
a_{_{\!1}},\ldots\!, a_{_{\!n}})D = {\boldsymbol c}(f)D 
		\;\;(\mbox{hence, } \Kr(D) \cap K =D )\,.
		\end{array}
		$$
		\end{enumerate} 
		
\smallskip

Kronecker's classical theory  led to two different 
major extensions:

$\bullet$ \; Beginning from   1936 \cite{Kr2},  W. Krull   
generalized the Kronecker function ring to the more 
general context of \sl integrally closed domains, \rm by introducing 
ideal systems associated to particular star operations: the a.b.  
(\it arithmetisch brauchbar\rm ) star 
operations.

$\bullet$ \; Beginning from  1956 \cite{Nagata:1956},  M. Nagata   
investigated, \sl  for an arbitrary integral domain 
$D$, \rm the domain
 $\left\{{f}/{g} \, \mid  \, f, g \in D[X] \; \; \mbox{and} \; \;
{\boldsymbol c}(g) = D \right\}$,
which coincides with $\Kr(D)$  if (and only if) 
$D$ is a Pr\"ufer domain \cite[Theorem 33.4]{Gilmer:1972}.  

\medskip

We recall these two major extensions of Kronecker's 
classical theory, but first we fix the general notation that we 
use in the sequel.

\bigskip
\noindent \bf General notation. \rm   Let $\,D\,$ be an integral 
domain
with quotient field $\,K$.  \  Let   
$\,\boldsymbol{\overline{F}}(D)\,$  
represent the set of all nonzero $\,D$--submodules of $\,K\,$ and 
  $\,\boldsymbol{F}(D)\,$   the nonzero fractionary 
ideals of
$\,D\,$ (i.e. $\,E \in \boldsymbol{\overline{F}}(D)\,$ such that $\,dE
\subseteq D\,$, \ for some nonzero element $\,d \in D\,$). Finally, 
let  
$\,\boldsymbol{f}(D)\,$  be the finitely generated
$\,D$-submodules of $\,K\,$.  Obviously:
$$ \boldsymbol{f}(D) \subseteq \boldsymbol{F}(D) \subseteq 
\boldsymbol{\overline{F}}(D) \,.
$$

\smallskip

One of the major difficulties for generalizing Kronecker's theory 
is 
that  Gauss Lemma  for the content of polynomials holds for 
Dedekind domains (or, more generally, for 
Pr\"ufer domains), but not in general:

\medskip 
\noindent \bf Gauss Lemma.  \sl  Let $f, g \in D[X]$, 
where $D$ is an integral domain. If $D$ is a Pr\"ufer domain, 
then: 
$$
{\boldsymbol c}(fg) = {\boldsymbol c}(f){\boldsymbol c}(g)\,,$$
and conversely. \rm Cf. \cite[Corollary 28.5]{Gilmer:1972}.
 \rm

\bigskip

For general integral domains, we always have the inclusion of ideals 
${\boldsymbol c}(fg) \subseteq {\boldsymbol c}(f) {\boldsymbol 
c}(g)$.  We also have the following result which is weaker than the 
Gauss Lemma but more widely applicable.

\smallskip

\noindent \bf Dedekind--Mertens Lemma.  \sl Let $D$ be an integral 
domain and  
$f, 
g \in D[X]$. Let $m := \deg(g)$. Then
$$
{\boldsymbol c}(f)^m{\boldsymbol c}(fg) = {\boldsymbol 
c}(f)^{m+1}{\boldsymbol 
c}(g)\,.$$
\rm Cf. \cite[Theorem 28.1]{Gilmer:1972}.
 \rm 

\medskip

In order to overcome this obstruction to generalizing the definition 
of Kronecker's function rings, Krull introduced 
multiplicative ideal systems having a  cancellation property which 
mirrors Gauss's Lemma. 
These ideal systems 
can  be defined by what are called now the  e.a.b.  (\it endlich 
arithmetisch brauchbar\rm ) star operations.

\medskip

\noindent \bf Star operations. \rm   A mapping   $\,\star : 
\boldsymbol{{F}}(D) \rightarrow
\boldsymbol{{F}}(D)\,$ ,  $\,I \mapsto I^\star\,,$\;
is called {\it a star operation of $\,D\,$} \rm 
if, for
all $\,z \in K\,$, $\,z \not = 0\,$ and for all $\,I, J \in
\boldsymbol{{F}}(D)\,$, the following properties hold:

 \hspace{10pt} $\mathbf{ \bf (\star_1)} \; \; \; (zD)^\star = 
zD\,, \;\;   (zI)^\star = 
zI^\star \,; $
  
 \hspace{10pt} $\mathbf{ \bf (\star_2)} \; \; \;   I \subseteq J 
\;\Rightarrow\; I^\star
\subseteq
J^\star \,;$

 \hspace{10pt} $\mathbf{ \bf (\star_3)} \; \; \;   I \subseteq 
I^\star 
 \; \; \textrm {   and  
}\; \; I^{\star \star} := (I^\star)^\star = I^\star\,. $

An \it e.a.b. star operation \rm on $D$  is a star operation $\star$ 
such that,  for all nonzero finitely 
generated ideals $I,\ J,\ H$ of $D$:
$$
(IJ)^\star \subseteq (IH)^\star \; \Rightarrow\; J^\star \subseteq 
H^\star\,.
$$
  \medskip
  
Using these notions, Krull recovers a useful identity for the 
contents 
of 
polynomials:

\medskip

\noindent \bf Gauss--Krull Lemma. \sl Let $\star$ be an e.a.b. star 
operation on an integral domain 
$D$ (this condition implies that $D$ is an integrally closed domain 
\cite[Corollary 32.8]{Gilmer:1972}) 
and let $f,\ g \in D[X]$ then:
$$
{\boldsymbol c}(fg)^\star = {\boldsymbol c}(f)^\star{\boldsymbol 
c}(g)^\star\,.$$
 \rm Cf. \cite[Lemma 32.6]{Gilmer:1972}.

\begin{reem} \rm Krull  introduced the concept of a star operation in 
his   first
\it Beitr\"age paper \rm in 1936 \cite{Kr2}. He used the notation  `` 
$^{\prime}$--Operation '' 
(``Strich--Operation'')   for his generic operation. [In this paper 
you can find the
terminology `` $^{\prime}$--Operation '' in footnote 13 and in the 
title of Section 6,
among other places.] 
     
      The notation  `` $\ast$--operation '' (``star--operation'')  
arises from
Section 26 of the original version of  Gilmer's ``Multiplicative 
Ideal Theory''
(1968) \cite{Gilmer:1968}.  Robert Gilmer gave us this explication  \ 
$\ll$ I believe the reason I switched from `` $^{\prime}$--Operation 
'' to `` $\ast$--operation '' was because \  `` $^{\prime}$ ''  
was not so generic at the time: \ $I'$ \ was frequently used as the 
notation
for the integral closure of an ideal $I$, just as $ D'$  was used to 
denote the
integral closure of the domain $D$. (Such notation was used, for 
example,
in both Nagata's Local Rings and in Zariski-Samuel's two volumes.) 
$\gg$  
     
Moreover,   Krull  only considered the concept of an
 ``arithmetisch brauchbar (a.b.) $^{\prime}$--Operation'', not of an 
e.a.b. operation.\\
Recall that an  \it a.b.--operation \rm  is a star operation $\star$ 
such 
that, if  $I \in \boldsymbol{{f}}(D)$ and $J, K \in 
\boldsymbol{{F}}(D)$ 
and  if $(IJ)^{\star} \subseteq (IH)^{\star}$ then $J^{\star} 
\subseteq 
H^{\star}$ .

The e.a.b. concept stems from the original version of
Gilmer's book \cite{Gilmer:1968}. The results of Section 26  show 
that this 
(presumably) weaker
concept is all that one needs to develop a complete theory of 
Kronecker
function rings. 
     
     In this regard, Robert Gilmer gave us this explication \  $\ll$ 
I believe I was influenced to recognize this because during the 1966 
calendar year in our graduate algebra seminar (Bill
Heinzer, Jimmy Arnold, and Jim Brewer, among others, were in that
seminar) we had covered Bourbaki's Chapitres 5 and 7 of 
\it Alg\`ebre Commutative\rm , and the development in Chapter 7 on 
the $v$--operation indicated
that e.a.b. would be sufficient.$\gg$   
\end{reem}

One of the main goals for the classical theory of star operations 
has been to construct a
Kronecker function
ring   associated to a domain, in a more general context than the 
original one
 considered by { L. Kronecker} in  1882.

More precisely, using star operations, in   1936  W.  Krull 
\cite{Kr2} defined a 
Kronecker function ring in a 
more general setting than Kronecker's.  (Further 
references are  H. Pr\"ufer (1932) \cite{Prufer:1932}, 
 Arnold  (1969) \cite{Arnold:1969},  Arnold-Brewer  (1971) 
\cite{Arnold-Brewer:1971},    
Dobbs-Fontana 
 (1986) \cite{Dobbs-Fontana:1986},  D.F. 
Anderson-Dobbs-Fontana  (1987) \cite{ADF:1987},
 Okabe-Matsuda   (1997) \cite{OM3}.)

\medskip

\noindent \bf Star--Kronecker function ring. \rm Let $D$ be an 
integrally closed integral domain with quotient 
field 
$K$ and let $\star$ be an 
e.a.b. star operation on $D$, then:
$$  
\Kr(D, \star) := \left\{\frac{f}{g} \, \mid \,   f, g \in D[X]\,, \; 
g \neq 0 \; \; 
\mbox{and} \; \; 
{\boldsymbol c}(f)^\star \subseteq {\boldsymbol 
c}(g)^\star\right\}$$ 
is an integral domain with quotient field $K(X)$, called \it  the 
$\star$--Kronecker function ring of $D$, having the following 
properties:
\begin{enumerate}
	\bf \item[(1)] \sl 
 $\Kr(D, \star)$ is a B\'ezout domain   and  $D[X] \subseteq \Kr(D, 
 \star) \subseteq 
	K(X)$.
	
		\bf \item[(2)] \sl Let $a_{_{\!0}}, a_{_{\!1}},\ldots, 
a_{_{\!n}}\in 
		D$ and set $f := a_{_{\!0}} +a_{_{\!1}}X+ \ldots + a_{_{\!n}}X^n 
\in 
		D[X]$, then:


		$$
		\begin{array}{rll}
		&(a_{_{\!0}}, a_{_{\!1}},\ldots, a_{_{\!n}})\Kr(D, \star) =&
		\hskip -30pt
		f\Kr(D, \star)\,,\\
		&(a_{_{\!0}}, a_{_{\!1}},\ldots, a_{_{\!n}})\Kr(D, \star) \cap K = 
&  \hskip -9 pt
		((a_{_{\!0}}, a_{_{\!1}},\ldots, a_{_{\!n}})D)^\star \\
		&(\mbox{i.e. } f\Kr(D, \star) \cap K = ({\boldsymbol 
c}(f))^\star)\,. \end{array}
		$$
In particular, $ \Kr(D, \star) \cap K = D^\star= D\,.$
\end{enumerate}  
  \rm For the proof cf. \cite[Theorem 32.7]{Gilmer:1972}.
\medskip

\bigskip

\noindent \bf Nagata's generalization of the Kronecker function 
ring. \rm  
The following  construction is possible for any integral domain $D$ 
and, even, for an arbitrary ring $D$.
 $$ \Na(D) := D(X) :=\left\{\frac{f}{g} \, \mid  \, f, g \in D[X] 
\; \; \mbox{and} \; \;
{\boldsymbol c}(g) = D \right\}\,,
$$
and this ring is called the \it Nagata ring of 
$D$. \rm 
This notion is essentially due to {Krull} (1943) 
\cite{Krull:1943}.  Then this ring was studied in {Nagata}'s book 
(1962) \cite[Section 6, page 17]{Nagata:1962}, using the notation 
$D(X)$, and in
{Samuel's Tata volume} (1964) \cite[page 27]{Samuel:1964} (where the 
notation $D(X)_{\mbox{\footnotesize loc}}$ was used). We introduced 
the notation $\Na(D)$ that is convenient for generalizations.
\medskip


In general, $\Na(D)$ is not a B\'ezout domain. It is not difficult 
to see that \cite[Theorem 33.4]{Gilmer:1972} and 
\cite[Theorem 8]{Anderson:1976} \newline
\begin{itemize} 

\item {\sl  $\Na(D)$ is a B\'ezout domain if and only if $D$ is 
a 
Pr\"ufer domain.}  
 
 \item {\sl $\Na(D)$ coincides with $\Kr(D)$ if and only if $D$ 
is a 
Pr\"ufer domain.} 

\item {\sl Every ideal of $\Na(D)$ is extended from $D$ if and only if $D$ 
is a 
Pr\"ufer domain.} \rm

\end{itemize}

\medskip

 The interest in Nagata's ring $D(X)$ is due to the fact that 
this 
ring of rational functions has some
strong ideal-theoretic  pro\-per\-ties that $D$ itself need not have, while  
maintaining a strict relation with the ideal structure of $D$. 

  \bf  (a) \sl The map 
$P \mapsto PD(X)$ establishes a 1-1 correspondence between the 
maximal ideals of $D$ and the maximal   ideals of 
$D(X)$\rm . 

  \bf  (b) \sl  For each ideal $\,I\,$ of $\,D\,$, \newline  
\centerline{$\, ID(X) \cap D = I\,$, \;
 \;  $\, D(X)/ID(X) \cong (D/ID)(X)\,$;}
\newline 
\centerline{$\,  I\,$ is finitely generated if and only if 
 $\,ID(X)\,$is finitely generated.} \rm 

Among the new
properties acquired by $D(X)$  are the following:

  \bf  (c) \sl  the
residue field at each maximal ideal of $\,D(X)\,$ is infinite;

  \bf  (d) \sl an ideal contained
in a finite union of ideals is contained in one of them;

  \bf  (e) \sl
each finitely generated locally principal ideal is principal 
(therefore \newline  
$\Pic(D(X)) = 0$).   \rm

The proofs of the previous results  can be found in Arnold (1969) 
\cite{Arnold:1969} and Gilmer's book  \cite[Proposition 33.1, 
Proposition 5.8]{Gilmer:1972}  (for \bf (a)\rm , \bf (b)\rm , and \bf 
(c) \rm which is a consequence of \bf (b)\rm ),  
Quartararo-Butts (1975) \cite{QB:1975}  (for \bf (d)\rm ) and D.D. 
Anderson (1977) \cite[Theorem 2]{Anderson:1977} (for \bf (e)\rm ). 

\bigskip

	\section{Basic facts on semistar operations}  \rm

In  1994, Okabe and 
Matsuda \cite{OM2} introduced the more
flexible notion of semistar operation $\,\star\,$ of an integral
domain $\,D\,,$ \ as a natural generalization of the notion of star
operation, allowing $\, D \neq D^\star\,$ (cf. also \cite{OM1}, 
\cite{MS}, and \cite{MSu}).  

\medskip

\noindent \bf Semistar operations. \rm  \rm  A mapping   $\,\star : 
\boldsymbol{\overline{F}}(D) \rightarrow
\boldsymbol{\overline{F}}(D)\,$,  $\,E \mapsto E^\star\,$
is called {\it a semistar operation of $\,D\,$} \rm 
if, for
all $\,z \in K\,$, $\,z \not = 0\,$ and for all $\,E,F \in
\boldsymbol{\overline{F}}(D)\,$, the following properties hold:

 \hspace{10pt} $\mathbf{ \bf (\star_1)} \; \; \; (zE)^\star = 
zE^\star \,; $
  
 \hspace{10pt} $\mathbf{ \bf (\star_2)} \; \; \;   E \subseteq F 
\;\Rightarrow\; E^\star
\subseteq
F^\star \,;$

 \hspace{10pt} $\mathbf{ \bf (\star_3)} \; \; \;   E \subseteq 
E^\star 
 \; \; \textrm {   and  
}\; \; E^{\star \star} := (E^\star)^\star = E^\star\,. $ \newline    
When $D^\star=D$, we say that $\star$ is     \emph{a (semi)star
operation} of $D$\rm , since, restricted to $\boldsymbol{F}(D)$ it 
is   \sl  a star
operation of $D$\,.  \rm

\medskip

For star operations,   the notion of  
$\star$--ideal leads to the definition of a canonically 
associated ideal system. 

For semistar operations, we need a more general notion,
that coincides with the notion of $\star$--ideal, when $\star$ is a
(semi)star operation.

$\bullet$ \;  \rm A nonzero (integral) ideal $I$ of
$D$ is a   \emph{quasi--$\star$--ideal} [respectively,\, {\it
$\,\star$--ideal} \rm]  if $I^\star \cap D = I$ [respectively,\, if 
$I^\star  = I$]. 

$\bullet$ \;  \ A {\it
quasi--$\star$--prime} \rm [respectively,\, {\it
$\,\star$--prime} \rm] of $D$ is a quasi--$\star$--ideal [respectively,\,
an integral $\,\star$--ideal] of $\,D\,$ which is also a prime.

$\bullet$ \;  \ A  \it quasi--$\star$--maximal  \rm 
[respecti\-vely,\, {\it
$\,\star$--maximal} \rm] of $\,D\,$ is a maximal element in the set of 
all
proper quasi--$\star$--ideals [respectively,\, integral 
$\star$--ideals]
of $\,D\,.$

We denote by\ Spec$^\star(D)$ \ [respectively,  Max$^\star(D)$, 
QSpec$^\star(D)$, QMax$^\star(D)$] the set of all $\star$--primes
[respectively,\, $\star$--maximals, \, quasi--$\star$--primes,\,
quasi--$\star$--maximals] of $\,D\,$.

For example, it is easy to see that,  \sl if $I^\star \neq D^\star$, 
then $I^\star \cap D$ is a quasi--$\star$--ideal that contains $I$  
\rm  (in particular, a $\star$--ideal is a quasi--$\star$--ideal).
\newline
Note that: 

-- \, when $\,D = D^\star\,$ the notions
 of quasi--$\star$--ideal and  $\,\star$--ideal coincide;

 -- \, 
$I^\star \neq D^\star$ is equivalent to $I^\star \cap D \neq D$.

\medskip

As in the classical star-operation setting, we associate to a \sl 
semistar \rm  ope\-ra\-tion $\,\star\,$ of $\,D\,$ a new semistar 
operation
$\,\star_f\,$ as follows.  \ If $\,E \in
\boldsymbol{\overline{F}}(D)\,$ we set:
$$
E^{\star_f} := \cup \{F^\star \;|\;\, F \subseteq E,\, F \in
\boldsymbol{f}(D)
\}\,.
$$
We call {$\,\star_f\,$ \it the 
semistar
operation of finite type of $D$ \rm associated to $\,\star\,$}.

$\bullet$\; If $\,\star =
\star_f\,$,\ we say that $\,\star\,$ is {\it a 
semistar ope\-ra\-tion of
finite type of $\,D\,$}.   Note that $\,\star_f \leq \star\,$ and
$\,(\star_f)_f = \star_f\,$,\ so $\,\star_f\,$ is 
of finite type on $\,D\,.$ 

\smallskip

  The following result is in \cite[Lemma 2.3]{FL3}.

\begin{leem}\label{le:1} \sl Let $\,\star\,$ be a non-trivial semistar
operation of finite type on $D$.  \ Then
\bara
 \bf \item \sl Each proper quasi--$\star$--ideal 
is
contained in a quasi--$\star$--maximal.
 \bf \item \sl Each
quasi--$\star$--maximal is a quasi--$\star$--prime.
  \bf
\item \sl Set
$$
\Pi^\star : = \{P \in \mbox{\rm Spec}(D)\, | \, P \not = 0 \,,  
P^\star \cap D \not = D\}\,.
$$
Then \,\rm QSpec$^\star(D) \subseteq \Pi^\star$\, \sl and the
set of maximal elements
of $\;\Pi^\star,$\ denoted by $\; {\Pi^\star_{\mbox{\rm \tiny
max}}} ,$ \ is nonempty and coincides with \rm QMax$^\star(D)$.
\eara \end{leem} \rm

For the sake of simplicity, when $\,\star = \star_f\,$, \ we will 
denote
simply by\ {${\calM}(\star) \,$}, the nonempty set
{$\,\Pi^\star_{\mbox{\rm \tiny max}} = 
\mbox{QMax}^\star(D)\,$}.

\bigskip
	
	\section{Nagata semistar domain}
	
	A generalization of the classical Nagata ring construction was
considered by  {Kang} (1987 \cite{Kang:1987} and 1989 
\cite{Kang:1989}). 
We further generalize this construction to the semistar setting.

\medskip

\bf Nagata semistar function ring. \rm  Given {{\sl any} integral 
domain} $\,D\,$ and 
{{\sl any} 
semistar operation $\,\star\,$} on
$\,D\,,$ \ we define {\it the semistar Nagata ring \rm 
} as
follows: 
$$
{\textstyle\rm Na}(D,\star) :=\left\{\frac{f}{g} \,  \mid \, f, g \in 
D[X] 
\,, \; g\neq 0\,,\; {\boldsymbol 
c}(g)^\star =D^\star \right\}\,.
$$

Note that $\,{\textstyle\rm Na}(D,\star) = 
{\textstyle\rm
Na}(D,\star_f)$\,.\  
Therefore, the assumption $\,\star = \star_f\,$ is not really 
restrictive when
considering Nagata semistar rings.
\medskip

If $\,\star = d\,$ is the identity
(semi)star operation of $D$, then:
$$
{\textstyle\rm Na}(D,d) = D(X)\,.
$$

	Some results on \sl star \rm Nagata rings proved by 
{Kang} in 1989 are
generalized to the semistar setting in the following:

\begin{prro}\label{prop:2}
\sl Let $\,\star\,$ be a nontrivial semi\-star
operation of an integral domain $\,D\,$. \ Set:\\
\centerline{
 $N(\star) := N_D(\star) := \{h \in D[X] \, |\; 
\boldsymbol{c}(h)^\star = 
 D^\star \}$ \,.}
$\mbox{  }$ \vskip -20pt  \bara
\rm \bf \item \sl $N(\star) = D[X] \setminus \cup \{Q[X] \; | \;\, Q 
\in
{\calM}(\star_f) \}\,$ is a saturated multiplicatively closed
subset of
$\,D[X]\,$ and $\, N(\star) = N(\star_f)\,$.

\rm \bf \item \rm  Max$(D[X]_{N(\star)}) = \{Q[X]_{N(\star)} \; |\;\, 
Q
\in
{\calM}(\star_f)\}\,.$

 \rm \bf \item ${\textstyle\rm Na}(D,\star) = D[X]_{N(\star)} = 
\cap\{ D_Q(X) \,| \, Q \in
 {\calM}(\star_f) \}\,.$

\rm \bf \item \sl ${\calM}(\star_f)$ coincides with the canonical  
image  in
\rm \, Spec$(D)$ \, \sl of the maximal  spectrum of
$\,{\textstyle\rm Na}(D,\star)\,;$ \,
  \sl i.e. 
${\calM}(\star_{f}) = \{ M \cap D \, | \, M \in $ \rm
Max$({\textstyle\rm Na}(D,\star)) \}\,. $

\eara
\end{prro}\rm 

For the proof cf. \cite[Theorem 3.1]{FL3}. \
From the previous Proposition \ref{prop:2} (4) we have:

\smallskip

\begin{coor} \label{cor:3}
	\sl Let $\,D\,$ be an integral domain, then:
	
\centerline{$Q \; \mbox{is a maximal $t$--ideal of } D \, 
\Leftrightarrow
\, Q = M \cap D\,, \; \mbox{for some} \, M \in \mbox{\rm Max(Na}(D, 
v))\,. $}

\end{coor} \rm

\bigskip


\section{The semistar operation associated to Na$(D, \star)$}

We start by recalling some distinguished classes of
semistar operations.

{$\bullet$} \, If $\,\Delta\,$ is a nonempty set of prime ideals of 
an integral 
domain
$\,D\,$, \ then the semistar operation {$\,\star_\Delta\,$}
defined on $\,D\,$ as follows
$$
E^{\star_\Delta} := \cap \{ED_P \;|\;\, P \in \Delta\}\,,
\; \; \textrm {  for each}    \; E \in \boldsymbol{\overline{F}}(D)\,,
$$
 \noindent is called {\it the spectral 
semistar
operation associated to \rm $\,\Delta\,$}.

{$\bullet$} \, A semistar operation 
$\,\star\,$ of an
integral domain $\,D\,$ is called {\it a spectral
semistar operation} \rm if there exists $\,\emptyset \not = \Delta
\subseteq \mbox{\rm Spec}(D)\,$ such that $\,\star = \star_\Delta\,$.

{$\bullet$}\, We say that $\,\star\,$ 
\it
{possesses enough primes} \rm or that $\,\star\,$ is { 
\it
a quasi-spectral semistar operation of $\,D\,$} \rm if, for each 
nonzero
ideal $\,I\,$ of $\,D\,$ such that $\,I^\star \cap D \not = D\,$, \
there exists a quasi--$\star$--prime $\,P\,$ of $\,D\,$ such that $\,I
\subseteq P\,$.

{$\bullet$}\, Finally, we say that 
$\,\star\,$ is
{\it a stable semistar operation on $\,D\,$} \rm if
$$
(E \cap F)^\star = E^\star \cap F^\star,  \;\, \textrm {  for all} \;
E,F \in \boldsymbol{\overline{F}}(D) \,.
$$

\begin{reem} \sl Mutatis mutandis \rm the previous notions were 
considered first in the star settings and, in particular, by  D.D. 
Anderson, D.F. Anderson and S. J. Cook  who gave important 
contributions to the subject \cite{Anderson:1988}, 
\cite{Anderson-Anderson:1990} and \cite{Anderson-Cook:2000}. The 
general situation was considered among the others by Fontana-Huckaba 
\cite{FH}, Fontana-Loper \cite{FL3}  and Halter-Koch \cite{HK4}.
\end{reem}


\begin{leem}\label{le:4} 
\sl Let $\,D\,$ be an integral domain and let
$\; \emptyset \not = \Delta
\subseteq \mbox{\rm Spec}(D)\,$. \ Then:

\bara
\rm \bf \item \sl $E^{\star_\Delta}D_P = ED_P$, \; for each $E \in
\boldsymbol{\overline{F}}(D)\,$ and for each $\,P \in \Delta\,$.  \rm
\bf \item \sl $(E \cap F)^{\star_\Delta}= E^{\star_\Delta} \cap
F^{\star_\Delta}$, \; for all $E,F \in \boldsymbol{\overline{F}}(D)$. 
\rm \bf \item \sl $P^{\star_\Delta} \cap D = P$, \, for each $P \in
\Delta$.  \rm \bf \item \sl If $\,I\,$ is a nonzero integral ideal of
$\,D\,$ and $\,I^{\star_\Delta} \cap D \not = D,\,$\,  then there exists
$\,P \in \Delta\,$ such that $\,I \subseteq P\,$.
\eara
 \end{leem} \rm
 \vskip -4pt For the proof cf. \cite[Lemma 4.1]{FH}.

\begin{leem}\label{le:5}
 \sl Let $\,\star\,$ be a nontrivial semistar o\-pe\-ration of
an integral domain $D$.  Then:
\bara
\rm \bf \item \sl $\,\star\,$ is spectral if and only if $\,\star\,$ 
is
quasi-spectral and stable.

\rm \bf \item \sl Assume that $\,\star = \star_f\,$.  \ Then 
$\,\star\,$
is quasi-spectral and $\,\calM(\star) \not = \emptyset\,$.\eara 
\end{leem} \rm  
\vskip -4pt For the proof cf. \cite[Theorem 4.12]{FH} and \cite[Lemma 
2.5]{FL3}.

\begin{thee}\label{the:6} \sl Let $\,\star\,$ be a nontrivial semistar
operation and let $E \in \boldsymbol{\overline{F}}(D)\,$.
 Set \vskip 0.1cm
\centerline{$ \tilde{\star} := (\star_f)_{sp} := 
\star_{{\calM}(\star_f)} \,.$}
[\ $\tilde{\star}$ is called {\it the spectral semistar
operation associated to $\,\star\,$}.]\, 
Then:
\bara
\rm \bf \item \sl $E^{\tilde{\star}} = \cap \{ED_Q \;|\;\, Q \in
{\calM}(\star_f) \}\,$\, [and  $E^{\star_f} = \cap \{E^{\star_f}D_Q 
\;|\;\, Q \in
{\calM}(\star_f) \}\,$].

\rm \bf \item \sl $\tilde{\star}  \leq
\star_f\,$.

\rm \bf \item \sl  $E{\textstyle\rm Na}(D,\star) = \cap\{ED_Q(X)\; 
|\;\,
Q \in
{\calM}(\star_f)\}\,$,\; thus:\\  $E{\textstyle\rm Na}(D,\star) \cap 
K = \cap
\{ED_Q \; | \;\, Q \in {\calM}(\star_f) \}\,$.

\rm \bf \item \sl $E^{\tilde{\star}} = E{\textstyle\rm Na}(D,\star)
\cap K$\,.
\eara
\end{thee} \rm

\vskip -5pt For the proof cf. \cite[Proposition 3.4]{FL3}.

\medskip 
Proposition~\ref{prop:2} (4) assures that, when a maximal
ideal of
\,Na$(D,\star)\,$ is contracted to $\,D\,$, \ the result is exactly a
prime
ideal in $\,{\calM}(\star_f)\,$. 
This result can be reversed. Moreover, the semistar operation 
$\tilde{\star}$ generates the same Nagata ring as $\star$. 

\begin{coor} \label{cor:7} \sl Let $\,\star,\ \star_{1},\ 
\star_{2}\,$ 
be semistar operations of
an integral domain $\,D\,$.  \ Then:
\bara
\rm \bf \item \sl $\, \Max(\Na(D,\star))  = \{QD_Q(X) \cap 
\Na(D,\star) \, \mid \, Q
\in {\calM}(\star_f)\}\,.$ 

\rm \bf \item \sl $(\tilde{\star})_f = \tilde{\star} =
 \tilde{\tilde{\star}}\,$.

\rm \bf \item \sl ${\calM}(\star_f) = {\calM}(\tilde{\star})\,$.

\rm \bf \item \sl $\Na(D,\star) = \Na(D,\tilde{\star}) \,$.

\rm \bf \item \sl \   $ \star_{1}\leq \star_{2} \; \Rightarrow \; 
\Na(D, \star_{1}) \subseteq  \Na(D, \star_{2}) \; \Leftrightarrow \; 
\widetilde{\ \star_{1}} \leq \widetilde{\ \star_{2}}$\,. 
\eara
\end{coor}
\vskip -4pt For the proof cf. \cite[Corollary 3.5 and Theorem 
3.8]{FL3}.

 \begin{reem} \label{rem:8} \rm Note that, when $\,\star\,$
is the (semi)star
{$v$--operation}, {then 
the
 (semi)\-star operation $\,\tilde{v}\,$ coincides with \it the 
(semi)star
 operation   $\,w\,$} \rm  defined as follows:  $$E^{w} :=
\cup \{(E:H)\; | \;\,H \in \boldsymbol{f}(D) \mbox{ and } H^v = D 
\}\,,$$
 $\mbox{for each } E \in \boldsymbol{\overline{F}}(D)\,,$ \ cf. 
\cite[page 182]{FH}.
     This (semi)star operation was   considered by 
{J. Hedstrom} and {E. Houston} in  1980  under 
the name of
     F$_{\infty}$--operation \cite{Hedstrom-Houston: 1980}.


	 Later, starting in  1997, this operation was
    studied by  {Wang Fanggui} and 
{R. McCasland} under the name of
    { $w$--operation} \cite{Fanggui-McCasland:1997} (cf. also 
   \cite {Fanggui:1}, \cite{Fanggui:2}
and \cite{Fanggui-McCasland:1999}).  (Unfortunately, the same 
notation is also used for the star a.b. operations defined by a 
family of valuation overrings \cite[page 398]{Gilmer:1972} and the 
two notions are not related, in general.)
    Note also that the notion of
     $w$--ideal coincides with the notion of  semi-divisorial 
     ideal   considered
     by {S. Glaz} and {W. 
Vasconcelos} in
      1977 \cite{Glaz-Vasconcelos:1977}.

	 Finally, in  2000, for
     each (semi)star operation $\,\star\,$, {D.D. 
Anderson} and {S.J. Cook} \cite{Anderson-Cook:2000}
      considered the {$\,\star_{w}$--operation} which
     can be defined as follows: 
	 $$E^{\star_{w}} := \cup \{(E:H)\; | \;\,H \in
     \boldsymbol{f}(D) \mbox{ and } H^\star = D \}\,,$$
	 $\mbox{for each } E \in \boldsymbol{\overline{F}}(D)\,.$ From
	 their theory (and from the results by Hedstrom and Houston) it 
follows that:
	 $$\,\star_{w} = \tilde{\star}\,.$$

	 The relation between
     $\,\tilde{\star}\,$ and {the localizing 
systems of ideals} (in the sense
     of Gabriel and Popescu) was established by {M. 
Fonta\-na}
     and {J. Huckaba} in  2000 \cite{FH}.
\end{reem}

	\section{The Kronecker function ring in a general setting} \rm

The problem of the construction of a Kronecker function ring for 
general integral domains was considered indipendently by  F. 
Halter-Koch  (2003)  \cite{HK3} and   Fontana-Loper  (2001, 2003) 
\cite{FL1}, \cite{FL3}.

Halter-Koch's approach is axiomatic and makes use of the theory of 
finitary ideal systems (star operations of finite type) \cite{HK2}. 
He also 
establishes a connection with  Krull's 
	theory of 
	Kronecker function rings and introduces the Kronecker function rings 
	for integral domains with an 
	ideal system which does not necessarily verify the cancellation 
	property (e.a.b.). \ Fontana-Loper's treatment is based on the 
Okabe-Matsuda's theory of semistar 
operations \cite{FL1}, \cite{FL3}, \cite{OM3}, and 
\cite{Matsuda:1998}.

 Halter-Koch  \cite{HK3} gives the following abstract definition 
 which does not rely on semistar operations or valuation overrings.

 \medskip
 
\bf $\boldsymbol K$--function ring. \rm  
 Let $K$ be a field, $R$  a subring of $K(X)$  and $D:= R \cap 
K$. If 

\bf (Kr.1) \rm  \ $X \in \boldsymbol{\mathcal{U}}(R)$ \, (i.e. $X$ is 
a unit in $R$\,);

\bf (Kr.2) \rm \ $f(0) \in fR$ for each $f \in K[X]$ \,;

\noindent then $R$ is called  \it a $K$--function ring of $D$.  \rm   

\medskip

Using only these two axioms, he proved that $R$ ``behaves as a 
Kronecker function ring'':

\begin{thee}  \sl Let $R$ be a $K$--function ring of $D=R\cap K$, 
then:
	\bara
	\rm \bf \item \sl $R$ is a B\'ezout domain with quotient field 
	$K(X)$\,.
	
	\rm \bf \item \sl  $D$ is integrally closed in $K$\,.
	
	\rm \bf \item \sl  For each polynomial  $\ f:= a_{0}+ a_{1}X +\ldots 
+a_{n}X^n 
\in 
	K[X]$,\ we have~  \newline \ $(a_{0}, a_{1}, \ldots, a_{n})R = fR$\,. 
	\eara
	
	\end{thee}
	
	 \vskip -4pt \rm For the proof cf. \cite[Theorem 2.2]{HK3}.
	 \medskip
	 
	 Our next goal is to describe Fontana-Loper's approach and to illustrate 
the relation with Halter-Koch's $K$--function rings.
	 
	 \medskip

  \bf Semistar Kronecker function ring. \rm  If {$\,\star\,$ is \sl 
any \rm semistar
operation} of { \sl any \rm  integral domain 
$\,D\,$},\
then we define {\it the Kronecker function ring 
of $\,D\,$
with respect to the semistar operation $\,\star\,$} \rm by:
$$
\begin {array} {rl}
\mbox{Kr}(D,\star) := \{ f/g  \, \ |& \, f,g \in D[X], \ g 
\neq 0, 
\;
\mbox{ and there exists }\\ & \hskip -4pt h \in D[X] \setminus \{0\} 
\; \mbox{ 
with } (\boldsymbol{c}(f)\boldsymbol{c}(h))^\star
\subseteq (\boldsymbol{c}(g)\boldsymbol{c}(h))^\star \,\}.
\end{array}
$$ 

At this point, we need some preliminaries in order:\\
--- to show that this construction leads to a natural extension 
of the classical Kronecker function ring,\\
--- to investigate the connections between the semistar Kronecker 
function ring  $\Kr(D,\star)$ and
the axiomatically defined $K$-function ring, \,  \\
--- to show that $\Kr(D,\star)$  defines a new semistar operation 
on $D$,  
behaving with respect $\Kr(D,\star)$ in a similar way to 
$\tilde{\star}$ with respect to $\Na(D, 
\star)$.

\bigskip
We start by recalling that it is possible to associate to an 
arbitrary semistar ope\-ration  an e.a.b. semistar operation.
\medskip

$\bullet$\; Given any semistar ope\-ration $\,\star\,$ of
$\,D\,$, \  we can define an e.a.b. semistar operation of finite type 
{$\,
\star_a\, $} of $\,D\,$, called {\it the e.a.b.
semistar
operation associated to $\,\star\,$}, \  as follows 
for each $ F \in \boldsymbol{f}(D)$ and
for each  $E \in {\overline{\boldsymbol{F}}}(D)\,:$
$$
\begin {array} {ll}
F^{\star_a} &:= \cup\{((FH)^\star:H^\star) \; \ | \; \, \; H \in
\boldsymbol{f}(D)\}\,, \\
E^{\star_a} &:= \cup\{F^{\star_a} \; | \; \, F \subseteq E\,,\; F \in
\boldsymbol{f}(D)\}\,. 
\end{array}
$$
 \rm  The previous
construction is essentially due to {P. Jaffard} 
(1960) \cite[Chapitre II, \S 2]{J}
and {F. Halter-Koch} (1997, 1998) \cite[Section 6]{HK1}, 
\cite[Chapter 19]{HK2} \ (cf. also Lorenzen (1939) 
\cite{Lorenzen:1939} and Aubert (1983) \cite{Aubert:1983}).

\smallskip

Obviously 
 $(\star_{_{\!f}})_{a}= \star_{a}$ . Note that (for instance 
\cite[Proposition 4.3 and 4.5]{FL1}): \\
--  \sl when $\star = 
\star_{_{\!f}}$, then $\star$ is e.a.b. if and only if $\star = 
\star_{a}\,.$ \\ \rm
--  {\sl $\,D^{\star_a}\,$ is integrally
closed and contains the integral closure of $\,D\,$}.\rm
 
 When $\,\star = v\,$, then {$\,D^{v_{a}}\,$} 
coincides
 with {\it the pseudo-integral closure of 
$\,D\,$ }\rm
 introduced by {D.F. Anderson}, 
{Houston}
 and {Zafrullah} (1992) \cite {AHZ}.

\begin{reem} \rm  In the classical context of \sl star \rm 
operations, $\star_a$ is expected to be a star operation too and for 
this reason is 
defined on the ``star 
closure'' of $D$ (or, on an integral domain which is ``star 
closed''), cf.   Okabe-Matsuda  (1992) \cite{OM1},   Halter-Koch 
 (1997, 1998, 2003) \cite{HK1}, \cite{HK2}, \cite{HK3}.\newline More 
precisely (even if $\star$ is a 
semistar operation), we call  \it the $\star$--closure of 
$D$\rm :
\\
\centerline{$ 
D^{\cl} := \cup \{(F^\star: F^\star) \mid F \in 
\boldsymbol{f}(D)\}\,.
$}\\


 It is easy to see that  \sl $D^{\cl}$  is an integrally closed 
overring 
of $D$ \rm and  $D$ is said  \it  $\star$--closed  \rm 
if 
 $D = 
D^{\cl}$. \newline
We can now define a new (semi)star operation 
on $D$ if $D =D^{\cl}$ (or, in general, a semistar operation on $D$), 
  cl$^\star$    by setting 
for each $F \in \boldsymbol{f}(D),$
for each  $E \in {\overline{\boldsymbol{F}}}(D)\,$:
$$
\begin {array} {rl}
F^{\cl} :=& \hskip -2pt \cup\{((H^\star:H^\star)F)^\star \; \ | \; \, 
\; H \in
\boldsymbol{f}(D)\}\,, \\
E^{\cl} :=& \hskip -2pt \cup\{F^{\cl} \; | \; \, F \subseteq E\,,\; F 
\in
\boldsymbol{f}(D)\}\,. 
\end{array} $$
 If we set $\overline{\star} :=\mbox{cl$^\star$}$, it 
is not 
difficult to see that  $D^{\cls} = D^{\cl}$  (and 
that it 
coincides with $D^{\star_{a}}$) 
 and $ D^{\cl}$ contains the ``classical'' integral closure of 
$D$. Moreover (as semistar operations on $D$):
$$  
\star_{_{\!f}} \leq \mbox{cl$^\star$} \leq \star_{a}\,,\;\;\; \;   
 (\star_{_{\!f}})_{a}= 
(\mbox{cl$^\star$})_{a}= (\star_{a})_{a}= \star_{a}\,.$$
 \end{reem}

 We now turn our attention to the valuation 
overrings.  The
notion that we recall next is due to {P. Jaffard} 
(1960) \cite[page 46]{J} (cf. also  Halter-Koch  (1997) 
\cite[Chapters 15 and 18]{HK2}).

\medskip 

$\bullet$\; For a domain
$\,D\,$ and a semistar operation $\,\star\,$ on $\,D\,$, \ we say 
that 
a
valuation overring $\,V\,$ of $\,D\,$ is {\it a
$\,\star$--valuation overring of $\,D\,$} \rm provided $\,F^\star
\subseteq FV\,,$ \ for each $\,F \in \boldsymbol{f}(D)\,$.\  \\ Note 
that,
by definition the $\,\star$--valuation overrings coincide with the
$\,\star_{f}$--valuation overrings.

\begin{prro}\label{prop:11} \sl Let $\,D\,$ be a domain and let 
$\,\star\,$
be a semistar operation on $\,D\,$. 
\bara
\bf \item \sl The $\,\star$--valuation overrings also coincide with 
the
$\,\star_{a}$--valuation overrings.
\bf \item \sl $D\cl = \cap \{V \mid  V \mbox{is a $\star$--valuation 
overring of } D\}\,.$
\bf \item \sl A valuation overring $\,V\,$ of
$\,D\,$ is a $\,\tilde{\star}$--valuation overring of $\,D\,$ if and
only if $\,V\,$ is an overring of $\,D_P\,,$ for some $\,P\in 
{\calM}(\star_f)\,$.
\eara
\end{prro} \rm 
\vskip -4pt For the proof cf. for instance \cite[Proposition 3.2, 3.3 
and Corollary 3.6]{FL2} and \cite[Theorem 3.9]{FL3}.

\begin{thee}\label{Theor:12} \sl Let $\,\star\,$ be a semi\-star 
operation of
an integral domain $\,D\,$ with quotient field $\,K\,.$ Then:
\bara
\rm \bf \item \sl ${\textstyle\rm Na}(D,\star) \subseteq 
{\textstyle\rm
Kr}(D,\star)\,.$

\rm \bf \item \sl $V$ is a $\star$--valuation overring of $D$ if and 
only if $V(X)$ is a valuation overring of $\Kr(D, \star)$.
\newline The map \ 
$W \mapsto W \cap K$ \ establishes a bijection between the set of all 
valuation overrings of $\Kr(D, \star)$ and the set of all the 
$\star$--valuation overrings of $D$.

\rm \bf \item \sl ${\textstyle\rm Kr}(D,\star) = {\textstyle\rm
Kr}(D,\star_f) = {\textstyle\rm Kr}(D,\star_a) = \cap \{V(X) \mid 
V$   
is a $\star$--valuation 
overring of $ D\}\,$ is a B\'ezout domain
with
quotient field $\,K(X)\,$.

\rm \bf  \item  \sl $E^{\star_a} = E{\Kr}(D,\star) \cap
K = \cap \{EV \mid  V   \mbox{is a $\star$--valuation 
overring of }$ $ D\}\,,$ for each $\,E \in 
\boldsymbol{\overline{F}}(D)\,.$ 

\rm \bf  \item  \sl $R:={\textstyle\rm Kr}(D,\star)$ is a 
$K$--function 
ring of $R\cap K = D^{\star_{a}}$ (in the sense of Halter-Koch's 
axiomatic 
definition).
\eara
\end{thee} 
\vskip -4pt For the proof cf. \cite[Theorem 3.11]{FL1}, \cite[Theorem 
3.5]{FL2}, \cite[Proposition 4.1]{FL3}.


\section{Some relations between Na$(D, \star)$\,,\; 
Kr$(D, \star)$\,,\;  and the semistar operations
	$\widetilde{\star}$\,,\;$\star_{a}$}  \rm

An elementary first question to ask is whether the two semistar
operations $\, \tilde{\star}\,$ and $\, \star_a \,$  are actually the 
same -
or usually the same - or rarely the same.

Proposition~\ref{prop:11} indicates that for a semistar operation
$\,\star\,$ on a domain $\,D\,$, the $\,\tilde{\star}$--valuation
overrings of $\,D\,$ are all the valuation overrings of the
localizations of $\,D\,$ at the primes in $\,{\calM}(\star_f)\,$.
On the other hand, we know from Theorem \ref{Theor:12} that the 
$\,\star_a$--valuation overrings (or,
equivalently, the $\,\star$--valuation overrings) of $\,D\,$ 
correspond
exactly to the valuation overrings of the Kronecker function ring
\,Kr$(D,\star)\,$. \  In particular, each $\,\star_a$--valuation 
overring is 
also a  $\,\tilde{\star}$--valuation
overring.

It is easy to imagine that these two collections of
valuation domains can frequently be different and, even when the two 
collections of valuation domains coincide, it may happen that 
$\,\tilde{\star} \neq
\star_{a}\,.$ \  Fontana-Loper  \cite{FL3} gives some examples which  
illustrate the different situations that can occur.

It is possible to prove positive statements about the relationship
between
$\,\widetilde{(\mbox{-})}\,$ and $\,(\mbox{-})_a\,$ under certain conditions
\cite[Proposition 5.4 and Remark 5.5]{FL3}.
However, we limit ourselves to stating a result 
that generalizes the fundamental result that is at the basis of 
Krull's theory of Kronecker function rings, i.e.
$ \Na(D) = \Na(D, d) = \Kr(D, b) = \Kr(D)$ if and only if $D $
is a Pr\"ufer domain, cf. for instance 
\cite[Theorem 33.4]{Gilmer:1972}.

\medskip

We recall the following definition, which generalizes the classical 
notion of Pr\"ufer domain.

\medskip

\bf Pr\"ufer semistar multiplication domain.  \rm Let $\star$ be a 
semistar operation on an integral domain $D$. A \it Pr\"ufer 
$\star$--multiplication domain 
 \rm (for short, a \it  P$\star$MD\rm ) \ is an integral domain 
$D$ such that $ (FF^{-1})^{\star_{_{\!f}}} = 
D^{\star_{_{\!f}}} \ (= D^\star)$\;  (i.e.,\  each $F$ is \it 
 ${\star_{_{\!f}}}$--invertible\rm ) for each $F \in \boldsymbol{f}(D)$

\medskip

Some of the statements of the following theo\-rem, due to 
Fontana-Jara-Santos (2003) \cite{Fontana-Jara-Santos:2003} 
gene\-ra\-li\-ze 
some of the classical
characterizations of  the Pr\"ufer $v$--multiplication domains  
(for 
short,  P$v$MD).

\begin{thee}\label{Theor:15} \sl Let 
$D$ be an integral
domain and $\star$ a semistar operation on $D$.  The following are
equivalent: \brom \rm \item \sl $D$ is a P$\star$MD.
\rm \item \sl
$\mbox{\rm Na}(D,\star)$ is a Pr\"{u}fer domain.
\rm \item \sl
$\mbox{\rm Na}(D,\star)=\mbox{\rm Kr}(D,{\star})$\,.
\rm \item \sl
$\tilde{\star} =\star_{a}$\,.
\rm \item \sl
${\star_{_{\!f}}} $ is stable and e.a.b..
\erom
\noindent In particular, $D$ is a P$\star$MD if and only if it is a
P$\tilde{\star}$MD. 
\end{thee}

\vskip -4pt For the proof cf. \cite[Theorem 3.1 and Remark 
3.1]{Fontana-Jara-Santos:2003}.

\medskip

The following gives the converse of the implication \  P$v$MD 
$\Rightarrow$ P$w$MD \  proved by  Wang-McCasland  
(1999) \cite[Section 2, page 160]{Fanggui-McCasland:1999}, cf. also  
D.D. Anderson-Cook (2000) \cite[Theorem 2.18]{Anderson-Cook:2000}.

\begin{coor}\label{Cor:16} \sl Let $D$ be an integral
domain. The following are
equivalent: \brom \rm \item \sl $D$ is a P$v$MD.
\rm \item \sl
$\mbox{\rm Na}(D, t)=\mbox{\rm Kr}(D, t)$\,.
\rm \item \sl
$w: = \tilde{v} =v_{a}$\,.
\rm \item \sl
$t $ is stable and e.a.b..
\erom
\noindent In particular, $D$ is a P$v$MD if and only if it is a
P$w$MD. 
\end{coor}

\vskip -4pt For the proof cf. \cite[Corollary 
3.1]{Fontana-Jara-Santos:2003}.

\medskip

In the star setting the relation between the P$\star$MDs and the 
P$v$MDs is described by the following:

\begin{coor}\label{cor:17} \sl Let $D$ be an integral
domain and $\star$ a  \rm star \sl  operation on $D$.
$$
D \mbox{ is a P$\star$MD } \, \Leftrightarrow \, D \mbox{ is a 
P$v$MD \ and \ } t =\widetilde{\star} \;  \mbox{ (or, equivalently, 
} t 
=\star_{_{\!f}} ).$$
\end{coor} \rm
\vskip -4pt For the proof cf. \cite[Proposition 
3.4]{Fontana-Jara-Santos:2003}.

\begin{reem}   \rm The  P$v$MDs were studied by Griffin in 1967 
\cite{Griffin:1967} under the name of $v$--multiplication domains, 
cf. also \cite{Kr1} and \cite{J}.
  Relevant contributions to the subject were given among the others 
by Arnold-Brewer  (1971) \cite{Arnold-Brewer:1971},  Heinzer-Ohm 
(1973) \cite{Heinzer-Ohm:1973}, Mott-Zafrullah  (1981) 
\cite{Mott-Zafrullah:1981}, Zafrullah  (1984) 
\cite{Zafrullah:1984},   Houston (1986) \cite{Houston:1986}, Kang  
(1989) \cite{Kang:1989},  Dobbs-Houston-Lucas-Zafrullah (1989) 
\cite{DHLZ:1989} and El Baghdadi (2002) \cite{ElBaghdadi:2002}.
  
  For $\star$ a star operation,  P$\star$MDs were considered 
by  
Houston-Malik-Mott in 1984  \cite{Houston-Malik-Mott:1984}, 
introducing a unified setting for studying Krull 
domains, Pr\"ufer domains and P$v$MDs.  This class of domains was 
also investigated by Garcia-Jara-Santos  (1999) 
\cite{Garcia-Jara-Santos:1999} and    
Halter-Koch  (2003) \cite{HK:2003} These papers led naturally 
to the study of the Pr\"ufer semistar multiplication domains 
initiated in 2003  by Fontana-Jara-Santos 
\cite{Fontana-Jara-Santos:2003}.

Related to this study are the questions on the invertibility property 
of ideals and modules especially in the star and semistar setting, 
cf. the survey paper by Zafrullah \cite{Zafrullah:2000}, Chang-Park 
(2003) \cite{Chang-Park:2003} and Fontana-Picozza (2005) \cite{Fontana-Picozza:2005}.

\end{reem}
 
\bigskip

\section{Intersections of local Nagata domains} \rm 

Given a semistar operation $\star$ on $D$, the 
	integral domains $\Na(D, \star)$  
	and  $\Kr(D, \star)$ (and the related semistar operations 
	\ $\tilde{\star}$ \ and \ $\star_{a}$) have in many regards a 
similar 
	behaviour. \ The following natural question is the starting point of 
a recent paper by M. Fontana and K.A. Loper \cite{FL4}: 
	
	\medskip   \sl
	Is it possible to find an integral domain of rational 
functions,   
 \rm denoted by \sl $\KN(D, \star)$ \ (\rm obtained  as  an 
intersection of local Nagata domains associated to   
any semistar operation $\star$\sl ), \ such that: \newline    
$\bullet$ \;  $ 
\Na(D, 
\star) \subseteq \KN(D, \star) \subseteq  \Kr(D, \star)$\,;  
\newline  
$\bullet$ \;  $\KN(D, \star)$  generalizes at the same time 
$\Na(D, 
\star)$ and $\Kr(D, \star)$ and coincides with 
$\Na(D, \star) = \Na(D, \tilde{\star})$ or $\Kr(D, \star)= \Kr(D, 
\star_{a})$, when the semistar 
operation (of finite type) $\star$ assumes 
the extreme values of the interval $\tilde{\star} \leq \star 
\leq \star_{a}$\ ?  \rm 

\medskip

In order to present the anwer to the previous question we need to 
settle some terminology:

\medskip

$\bullet$\, If $F$ is in $\boldsymbol{f}(D)$, we say that $F$ is  \it 
             $\star$--e.a.b. \rm [respectively, \it   $\star$--a.b. 
\rm ]
if 
$(FG)^{\star}
             \subseteq (FH)^{\star}$ ,  with $G,\ H \in 
\boldsymbol{f}(D)$ [respectively, $G,\ H \in 
\overline{\boldsymbol{F}}(D)$],
 implies
             that $G^{\star}\subseteq H^{\star}$\rm .

\begin{leem} \label{le:18}  \sl Let $\star$ be a semistar operation 
on an integral domain
$D$, let $F \in \boldsymbol{f}(D) $ be ${\star_{_{\!f}}}$--invertible 
and let $(L,
N)$ be a local $\star$--overring of $D$.  Then $F\! L$ is a principal
fractional ideal of $L$.
\end{leem}  \rm
\vskip -4pt For the proof cf. \cite[Corollary 4.3]{FL4}.
\medskip

Note that, in general, $\star$--(e.)a.b. does not imply 
$\star$--invertible,
even for finite type semistar operations.  However, it is 
possible to 
show that,  for finite type stable semistar operations $\star$ (i.e. 
when 
$\star = \tilde{\star}$), the notions of  \ $\star$--e.a.b., \ 
$\star$--a.b. \  
and \ $\star$--invertible coincide  \cite[Proposition 5.3(2)]{FL4}.

 \medskip

Our next goal is to generalize Lemma \ref{le:18} to the case of 
$\star$--e.a.b. ideals.

\medskip

\bf Semistar monolocalities. \rm  Let 
         $\star$ be a semistar operation on an integral domain  
         $D$.
	  \it A $\star$--monolocality  
	      of $D$  \rm is  a local overring $L$  of $D$ such 
that: \newline
	      --   $FL$ is a principal fractionary ideal of $L$, \ 
	     for each $\star$--e.a.b. $F \in \boldsymbol{{f}}(D)$\,;  
	     \newline \medskip
	      --  $L = L^{\star_{_{\!f}}}$.

	   \medskip
	      
	      Obviously,  each $\star$--valuation overring is a 
	      $\star$--monolocality.  \
	      It is not hard to prove that,  for each $Q \in 
	      \calM({\star_{_{\!f}}})$, $D_{Q}$ is a 
	      $\tilde{\star}$--monolocality    \cite[Proposition 5.3(1)]{FL4}. \ Set:

$$\begin{array}{rl}  \calL(D, \star):=& \{ L\mid L \mbox{ is a } \star\mbox 
{--monolocality of } D \} \\
\KN(D, \star) := & \cap \{ L(X) \mid L \in 
\calL(D, \star)\}\,. \end{array}$$

	     We are now in condition to state some among the main results 
proved in Fontana-Loper (2005)   \cite{FL4}.	     

\begin{thee} \label{th:19} \sl Let $\star$ be a semistar operation 
on an
integral domain $D$ with quotient field $K$.  \bara
\bf \item \sl $ \Na(D, \star) \subseteq  \KN(D, 
\star)  \subseteq  \Kr(D, \star))\,.$

 \bf \item \sl $\begin{array}{rl}
    \KN(D, \star) :=  \{ f/g\in K(X) \mid  &  f, g\in
D[X],\ 
  g \neq 0,\ \mbox{ such that }  \boldsymbol{c}(g)\\ &\hskip -4pt  
\mbox{ is  }\star\mbox{--e.a.b.  and  } \boldsymbol{c}(f) \subseteq 
\boldsymbol{c}(g)^\star \}.
 \end{array} $
 
        \bf \item \sl For each maximal ideal 
${\boldsymbol{\mathfrak{m}}}$ of $\KN(D, \star)$, set
       $L({\boldsymbol{\mathfrak{m}}}) := 
       \KN(D, \star)_{{\boldsymbol{\mathfrak{m}}}} 
\cap K$.\  
Then: 
\smallskip \newline 
     $L({\boldsymbol{\mathfrak{m}}})$ is a $\star$--monolocality of
     $D$ (with maximal ideal\ $\boldsymbol{\mathscr{M}}:= 
{\boldsymbol{\mathfrak{m}}}  \KN(D,
     \star)_{{\boldsymbol{\mathfrak{m}}}} \cap 
     L({\boldsymbol{\mathfrak{m}}})$),
     \medskip \newline $\KN(D, 
     \star)_{{\boldsymbol{\mathfrak{m}}}}$
     coincides with the Nagata ring $L({\boldsymbol{\mathfrak{m}}}) 
(X)$  
     and 
     ${\boldsymbol{\mathfrak{m}}}$ coincides with 
$\boldsymbol{\mathscr{M}}(X) 
\cap \KN(D, \star)$.
     \bf \item \sl
   Every $\star$--monolocality of an
    integral domain $D$ contains a minimal $\star$--mono\-locality of
    $D$. If we denote by ${\calL}(D, \star)_{min}$  the set of all 
the minimal 
       $\star$--monolocalities of $D$, then\newline  
       ${\calL}(D, \star)_{min} =
       \{ L({\boldsymbol{\mathfrak{m}}}) \mid  
       {\boldsymbol{\mathfrak{m}}} \in \Max(\KN(D, 
\star)) \}$ \ and, obviously, \newline
       $ {\KN}(D, \star) = \cap 
\{L(X) \mid L 
    \in {\calL}(D, \star)_{min} \}\,.$
     \bf \item \sl For each $J:= (a_{0},\ a_{1},\ \ldots,\  a_{n})D 
\in
     \boldsymbol{f}(D)$, with $J\subseteq D$ and $J$ 
$\star$--e.a.b.,  let
     $g:=a_{0}+ a_{1}X+\ \ldots\ + a_{n}X^n \in D[X]$, then:
     $$
     J\KN(D, \star) =J^\star\KN(D, 
\star)=g\KN(D, \star)\,.$$
   
     \bf \item \sl 
   Let  $ 
\star_{_{\!\ell}}$ 
   and   ${\wedge_{\calL}} $ 
  be  the  semistar operations of $D$ 
defined as follows, for each $E \in \overline{\boldsymbol{F}}(D) $,
  $$\begin{array}{rl}  
   E^{\star_{_{\!\ell}}}:=&  E\KN(D, \star) \cap 
K\,,\\ 
     E^{\wedge_{\calL}}:=&   \cap \{EL \mid L\in 
{\calL}(D, 
    \star)\}\,.
    \end{array}
    $$ 
   Then: \newline \centerline{ $\widetilde{\star}  \leq 
{\star_{_{\!{\ell}}}} = \wedge_{\calL}   
 \leq \star_a \,.$}
     \bf \item \sl \; $
             \Na(D, \star) = \Na(D, \widetilde{\star}) = \KN(D, 
             \widetilde{\star})$\ and\ $\KN(D,  \star_{a}) =\Kr(D, 
\star_{a}) = \Kr(D, \star)\,.$
\eara
 \end{thee}
\vskip -4pt For the proof cf.  \cite[Theorem 5.11 ((1),(4),(5),(6), and (7)), Proposition 6.3, and 
Corollary 6.4 ((1) and (2)) ]{FL4}.

     \end{document}